\documentclass[a4paper,12pt]{article}
\usepackage{amsmath,amssymb}

\newtheorem{thm}{Theorem}
\newtheorem{lem}{Lemma}
\newtheorem{prop}{Proposition}
\newtheorem{cor}{Corollary}
\newtheorem{defn}{Definition}

\begin{document}
\title{On quasi-separative semigroups}
\author{Yu.\,I.\,Krasilnikova, B.\,V.\,Novikov}
\date{}
\maketitle

\texttt{ABSTRACT.}  We built some congruences on semigroups, from
where a decomposition of quasi-separative semigroups was obtained.

\section{Introduction}

The research of separative semigroups was being begun from the
famous paper of Hewitt and Zuckerman \cite{cl-pr,h-z}, where, in
particular, they proved that any commutative separative semigroup
is isomorphic to a semilattice of cancellative semigroups. An
generalization for the noncommutative case has been made by
Burmistrovich \cite{burm} and independently by Petrich
\cite{petr}. Drazin \cite{dr} introduced the term
`quasi-separativity' and studied connections between it and others
semigroup properties (inversity, regularity etc).

 We shall follow the terminology proposed by Drazin:
\begin{defn}
A semigroup $S$ is called separative\footnote{Burmistrovich
\cite{burm} called it {\it weakly cancellative}} if
\begin{equation}\label{eq1}
 \left\{
\begin{array}{llll}
     x^2=xy \\
     y^2=yx \\
     \end{array}
\right.\ \Longrightarrow\ x=y \qquad \mbox{\rm and} \qquad
 \left\{
\begin{array}{llll}
     x^2=yx \\
     y^2=xy \\
     \end{array}
\right.\ \Longrightarrow\ x=y
\end{equation}
for all $x,y \in S$.

A semigroup $S$ is called quasi-separative if
\begin{equation}\label{eq2}
x^2=xy=yx=y^2 \Longrightarrow x=y
\end{equation}
for all $x,y \in S$.
\end{defn}

Drazin also showed that in the definition of quasi-separativity we
can replace (\ref{eq2}) by the next condition
 \begin{equation}\label{eq3}
x^2=xy=y^2 \Longrightarrow x=y.
\end{equation}
It often simplifies considerably proofs of assertions.

The main result of this paragraph is an extension of the
Burmistrovich's theorem (Theorem \ref{10}): any quasi-separative
semigroup is decomposable into a semilattice of subsemigroups,
which are called {\it quasi-cancellative} by us. With this aim we
previously build certain congruences on arbitrary semigroup
(Theorem \ref{4}); they give semilattice decompositions in the
quasi-separative case. As a corrolary, in Sect.\,4 we consider an
intermediate class of semigroups ({\it weakly balanced
semigroups}) between separative and quasi-separative ones and
discuss the connections between them.

\section{Relation $\Omega$}

Let $S$ be an arbitrary semigroup. By analogy with [2], we define
two binary relations $E(a), F(a)\subset S\times S$ for every
element $a\in S$:
$$
E(a)=\{(x,y)\mid ax=ay\},\quad F(a)=\{(x,y)\mid xa=ya\}.
$$

The next properties of these relations are obvious:
\begin{eqnarray}
&\  E(b)\subset E(ab)\label{eq4}\\
&\  F(a)\subset F(ab)\label{eq5}\\
&\  bE(ab)\subset E(a)\label{eq6}\\
&\ F(ab)a\subset F(b)\label{eq7}
\end{eqnarray}
(here and below for a binary relation $R\subset S\times S$ and for
an element $x\in S$ the relation $\{(xa,xb)\mid (a,b)\in R\}$ is
denoted by $xR$; analogously, $Rx$).

In what follows, the main tools for our studying will be the
relations $\Omega\subset S\times S$, which satisfy the next
conditions:
\begin{eqnarray}
\forall a\quad\Omega\cap E(a)=\Omega\cap F(a)\label{eq8}\\
\forall a,b \quad b(\Omega \cap E(ab)) \subset \Omega \label{eq9}\\
\forall a,b \quad (\Omega \cap F(ab))a \subset \Omega \label{eq10}
\end{eqnarray}
and the equivalences $\sim_{\Omega}$ on $S$ corresponding to these
relations:
$$
a\sim_{\Omega} b\Longleftrightarrow \Omega\cap E(a)=\Omega\cap
E(b).
$$

According to (\ref{eq8}), such definition is equal to the
following:
$$
a\sim_{\Omega} b\Longleftrightarrow \Omega\cap F(a)=\Omega\cap F(b).
$$

\begin{lem}\label{3}
For all elements $a,b\in S$ and a relation $\Omega$ which
satisfies the conditions (\ref{eq8})-(\ref{eq10})
$$
\Omega\cap E(a)\subset\Omega\cap E(ab)\cap E(ba).
$$
\end{lem}
{\bf Proof.} By (\ref{eq8}) we have:
$$
\Omega\cap E(a)=(\Omega\cap E(a))\cap(\Omega\cap F(a)).
$$

From here, using (\ref{eq4}) and (\ref{eq5}), we get:
$$
\Omega\cap E(a)\subset (\Omega\cap E(ba))\cap (\Omega\cap F(ab))=\Omega\cap E(ba)\cap E(ab)
$$
(the last equality follows from (\ref{eq8})). $\blacksquare$

Our first result is fulfilled for an arbitrary semigroup:
\begin{thm}\label{4}
Let $\Omega\subset S\times S$ satisfies the conditions
(\ref{eq8})-(\ref{eq10}). Then the equivalence $\sim_{\Omega}$ is
a congruence on $S$.
\end{thm}
{\bf Proof.} Let $a, b, c\in S$ and $a\sim_{\Omega} b$. Obviously,
for the proving the right compatibility of $\sim_{\Omega}$ it is
enough to verify the inclusion
$$
\Omega\cap E(ac)\subset \Omega\cap E(bc).
$$

Let $(x, y)\in \Omega\cap E(ac)$. Owing to (\ref{eq9}) $(cx,cy)\in
\Omega$. On the other hand, (\ref{eq6}) implies the inclusion
$(cx, cy)\in cE(ac)\subset E(a)$. Therefore,
$$
(cx, cy)\in \Omega\cap E(a)=\Omega\cap E(b).
$$
From here it follows that $(x, y)\in \Omega\cap E(bc)$.

Similarly, by (\ref{eq7}) and (\ref{eq10}) the left compatibility
can be proved. $\blacksquare$

\smallskip

\noindent{\bf Example.} Let $S$ be a commutative semigroup,
$\Omega = S\times S$. Then the conditions of Theorem \ref{4} are
true and the equivalence
$$
a\sim b\Longleftrightarrow  E(a)= E(b)
$$
is a congruence relation.

\section{A decomposition of quasi-separative semi\-gro\-ups}

   In this section we apply the preceding theorem to
quasi-separative semigroups.

Note that the definition of quasi-separativity in the form
(\ref{eq3}) may be formulated in terms of the relations $E(a)$ and
$F(a)$:
\begin{equation}\label{eq11}
(a,b)\in E(a)\cap F(b) \Longrightarrow a=b
\end{equation}
for all $a, b\in S.$

\begin{thm}\label{6}
Let $\Omega$ be a relation on quasi-separative semigroup $S$ which
satisfies the conditions (\ref{eq8})-(\ref{eq10}). Then
$S/\!\!\sim_{\Omega}$ is a semilattice.
\end{thm}
{\bf Proof.} First, show that $S/\!\!\sim_{\Omega}$ is a band. In
order to verify this statement it is sufficient to justify that
the equality $\Omega\cap E(a)=\Omega\cap E(a^2)$ is right for any
$a\in S.$

An inclusion
$$
\Omega\cap E(a)\subset \Omega\cap E(a^2)
$$
at once follows from Lemma \ref{3}. Conversely, if $(x, y)\in
\Omega\cap E(a^2)$, then $(ax, ay)\in E(a)$. Moreover, owing to
(\ref{eq9})
$$
(ax, ay)\in a(\Omega\cap E(a^2))\subset \Omega,
$$
hence, $(ax, ay)\in \Omega\cap E(a)$. By Lemma \ref{3}
$$
(ax, ay)\in \Omega\cap E(ax)\cap E(xa)\cap E(ya)\cap E(ay),
$$
whence, in particular,
$$(ax, ay)\in \Omega\cap E(ax)\cap E(ay)=\Omega\cap E(ax)\cap F(ay)$$
by the condition (\ref{eq8}). From (\ref{eq11}) we obtain $ax=ay$,
that is $(x, y)\in E(a).$ Therefore, $\Omega\cap E(a^2)\subset
\Omega\cap E(a)$ and the first part of Theorem is proved.

Now we shall prove that $S/\!\!\sim_{\Omega}$ is commutative, viz.
that $\Omega\cap E(ab)=\Omega\cap E(ba)$. The successive using
the properties (\ref{eq4}), (\ref{eq8}) and (\ref{eq5}) gives us:
$$
\Omega\cap E(ab)\subset \Omega\cap E(bab)=\Omega\cap F(bab)\subset
\Omega\cap F((ba)^2)=\Omega\cap E((ba)^2).
$$

Since, as proved above, $S/\!\!\sim$ is a band, then
$$
\Omega\cap E(ab)\subset \Omega\cap E(ba).
$$

Analogously,
$$
\Omega\cap E(ba)\subset \Omega\cap E(ab),
$$
what completes the proof of Theorem. $\blacksquare$

Next assertion gives us a preliminary information about
$\sim_{\Omega}$-classes. Denote by $\Delta_T$ the diagonal of
Cartesian square $T\times T$.
\begin{prop}\label{7}
If $S$ is a quasi-separative semigroup, then each
$\sim_{\Omega}$-class $T$ satisfies the next condition for all
$a\in T$:
$$
\Omega\cap E(a)\cap (T\times T)\subset \Delta_{T}.
$$
\end{prop}
{\bf Proof.} Indeed, let $(x, y)\in \Omega \cap E(a)\cap (T\times
T).$ Since $x\sim_{\Omega} y\sim_{\Omega} a$, then
$$
(x,y)\in \Omega\cap E(a)=\Omega\cap E(x)\cap F(y),
$$
whence, by (\ref{eq11}), we have $x=y.$ $\blacksquare$

\begin{defn}\label{9}
We call a semigroup $S$ quasi-cancellative if the condition
$$
\left\{
\begin{array}{ll}
\forall x,y\in S^1 \quad xby=xcy \Longleftrightarrow yxb=yxc
\Longleftrightarrow byx=cyx \\
ab=ac.
\end{array}
\right.
$$
implies $b=c$.
\end{defn}

Obviously, every right- or left-cancellative semigroup is
quasi-cancellative.

Our main result on structure of quasi-separative semigroups is the
next
 \begin{thm}\label{10}
A semigroup is quasi-separative if and only if it is a semilattice
of quasi-separati\-ve quasi-cancellative semigroups.
\end{thm}
{\bf Proof.} {\it Necessity.} Denote a binary relation $\Omega_S$
on $S$:
\begin{equation}\label{eq12}
\Omega_S=\{(x,y)\mid \forall a, b\in S^1  \quad
axb=ayb\Longleftrightarrow
 xba=yba\Longleftrightarrow bax=bay\}
\end{equation}
and verify the conditions (\ref{eq8})-(\ref{eq10}) for it. Since
for $a=1$ we have:
$$
xb=yb\Longleftrightarrow bx=by,
$$
for any pair $(x, y)\in \Omega_S$, then obviously, (\ref{eq8})
holds. Now we prove that $\Omega_S$ is left compatibility, from
where (\ref{eq9}) will follow.

Let $(x,y)\in \Omega_S$, $b\in S$. To prove that $(bx,by)\in
\Omega_S$ one needs to check the fulfilment of the implications:
$$
\forall c,d\in S^1 \quad cbxd=cbyd\Longleftrightarrow
dcbx=dcby\Longleftrightarrow bxdc=bydc.
$$

The implication $cbxd=cbyd\Longleftrightarrow dcbx=dcby$
immediately follows from the definition of $\Omega_S$. Let
$dcbx=dcby$. From (\ref{eq12}) we obtain:
$$
dcbx=dcby\Longrightarrow xdcb=ydcb.
$$
Therefore
$$
(bxdc)^2=bx(dcbx)dc=b(xdcb)ydc=(bydc)^2
$$
and quasi-separativity implies $bxdc=bydc$.

Similarly, if $bxdc=bydc$, then
$$
bxdc=bydc\Longrightarrow dcbx=dcby \Longrightarrow xdcb=ydcb.
$$
Hence
$$
(cbxd)^2=cbx(dcbx)d=c(bxdc)byd=(cbyd)^2
$$
and $cbxd=cbyd$.

In the same way right compatibility is checked, and so the
condition (\ref{eq10}) is fulfilled.

Thus, $S/\!\!\sim_{\Omega_S}$ is a commutative band by Theorem
\ref{6}. It remains to show that its components are
quasi-cancellative.

Let suppose that the conditions of Definition \ref{9} hold for
some elements $a,b,c,d$ from the $\sim_{\Omega_S}$-class $T\subset
S$. It means that
$$
(c,d)\in \Omega_T\cap E(a)\subset \Omega_S\cap E(a)= \Omega_S\cap
E(b).
$$
Hence $bc=bd$. Moreover, (\ref{eq8}) implies $cb=db$. In
particular, replacing $b$ in the obtained equations by $c$ and
$d$, we get:
$$
c^2=cd=dc=d^2,
$$
whence $c=d$.

{\it Sufficiency.} It is easy to see that any semilattice of
quasi-separative semigroups is also quasi-separative.
$\blacksquare$

\section{Corollaries and Examples}

In this section we show that Theorem \ref{10} implies the theorem
of Burmistrovich on the separative semigroups and obtain an
assertion about certain intermediate class of semigroups.

\begin{prop}\label{11}
Every separative quasi-cancellative semigroup is cancellative.
\end{prop}
{\bf Proof.} Let $S$ be separative and quasi-cancellative,
$a,b,c\in S$, $ab=ac$. By Lemma 1 from \cite{burm} for all $x,y\in
S$
$$
xby=xcy\Longrightarrow byx=cyx \Longrightarrow yxb=yxc
\Longrightarrow xby=xcy.
$$
So, by quasi-separativity $b=c$.

To prove the right cancellativity we ought to apply the Lemma 1
\cite{burm} to the equality $ba=ca$ and to refer to the previous
argumentation. $\blacksquare$

\begin{cor}\label{12}{\bf (Burmistrovich's Theorem \cite{burm})}
A semigroup is separative if and only if it is isomorphic to a
semilattice of cancellative semigroups.$\blacksquare$
\end{cor}

\begin{defn}\label{13}
A semigroup $S$ is called \it{weakly cancellative} \cite{petr} if
for every $a,b,x,y \in S$
\begin{equation}\label{eq13}
 \left\{
\begin{array}{llll}
     ax=ay \\
     xb=yb \\
     \end{array}
\right.
\end{equation}
implies $x=y$.

We call a semigroup $S$ weakly balanced, if (\ref{eq13}) implies
$$
 \left\{
\begin{array}{llll}
     xa=ya \\
     bx=by. \\
     \end{array}
\right.
$$
\end{defn}

Obviously, every weakly cancellative semigroup is
quasi-separative; but in general this is not hold in the weakly
balanced case (for example, all commutative semigroups are weakly
balanced). On the other hand, by above-mentioned Lemma 1
\cite{burm} all separative semigroups are weakly balanced, so two
next facts give a partial extension of Burmistrovich's theorem to
the more wide class of semigroups.

\begin{prop}\label{14}
If $S$ is a quasi-cancellative weakly balanced semigroup, then $S$
is also weakly cancellative.
\end{prop}
{\bf Proof.} Let $S$ be quasi-cancellative and weakly balanced,
$a,b,x,y\in S$ and
$$
ax=ay,\qquad xb=yb.
$$
If $uxv=uyv$ for some elements $u,v\in S^1$, then by the weakly
balancity from this last equality and from $axv=ayv$ we obtain
$xvu=yvu$. Similarly, implications
$$
xvu=yvu\Longrightarrow vux=vuy \Longrightarrow uxv=uyv.
$$
can be obtained. Now $x=y$ because of quasi-cancellativity.
$\blacksquare$

\begin{cor}\label{15}
Every quasi-separative weakly balanced semigroup is isomorphic to
a semilattice of weakly cancellative semigroups. $\blacksquare$
\end{cor}

We don't know whether the converse to the Corollary \ref{15} is
true. One can affirm only that a semilattice of weakly
cancellative semigroups (which, evidently, is quasi-separative)
satisfies the next condition:
\begin{equation}\label{eq14}
\left\{
\begin{array}{cccc}
a^2x&=&a^2y \\
xa^2&=&ya^2
\end{array}
\right. \Longrightarrow \left\{
\begin{array}{cccc}
ax&=&ay \\
xa&=&ya
\end{array}
\right.
\end{equation}

Really, it follows out of $a^2x=a^2y$ that $ax,ay,xa,ya$ contain
in the same component of the semilattice. From the antecedent of
(\ref{eq14}) we have:
$$
\left\{
\begin{array}{llll}
     (xa)(ax)&=&(xa)(ay)\\
     (ax)(a^2x)&=&(ay)(a^2x) \\
     \end{array}%
\right.
$$
Now weakly cancellativity implies $ax=ay$ and, similarly, $xa=ya$.

\bigskip

In conclusion we discuss the connections between considered
classes of semigroups. They may be presented by a diagram:

$$
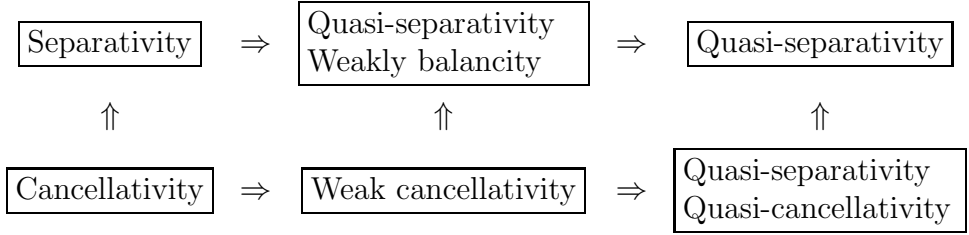
\begin{array}{ccccc}
\fbox{Separativity}&\Rightarrow &\fbox{\begin{minipage}{20ex}
Quasi-separativity
\\ Weakly balancity
\end{minipage}}
&\Rightarrow
&\fbox{Quasi-separativity}\\
\\[-2ex]
\Uparrow &&\Uparrow &&\Uparrow \\[1ex]
\fbox{Cancellativity}&\Rightarrow &\fbox{Weak cancellativity}
&\Rightarrow &\fbox{\begin{minipage}{20.1ex}Quasi-separativity\\
Quasi-cancellativity
\end{minipage}}
\end{array}
$$

\smallskip

Now we shall show that all implications in this picture are
strict.

Obviously, any commutative quasi-cancellative semigroup is
cancellative. Hence not every separative semigroup is
quasi-cancellative. From here it follows that all the vertical
implications are strict.

Every completely simple semigroup is weakly cancellative, but not
separative (if it is not a group). Hence the left horizontal
implications are strict.

 Bicyclic semigroup $B=\langle a,b\mid ba=1\rangle$ is
 quasi-separative. Since $B$ is simple, it cannot be
 decomposed into a nontrivial semilattice of its subsemigroups. By Theorem
\ref{10} it is quasi-cancellative. On the other hand, the
equalities
$$
b^2\cdot 1= b^2\cdot ab, \qquad 1\cdot a=ab\cdot a
$$
imply that $B$ is not weakly balanced. From here it follows that
the right horizontal implications are strict.

\bigskip

\begin{quote}
Lugansk National Pedagogical University

kji@leasat.net

\bigskip

Kharkov National University

boris.v.novikov@univer.kharkov.ua
\end{quote}

\end{document}